\documentclass[12pt,letterpaper]{amsart}
\usepackage{amssymb,latexsym, amsmath, amsxtra}
\usepackage[dvips]{graphics}

\theoremstyle{plain}
\newtheorem{theorem}{Theorem}[section]

\newtheorem{lemma}[theorem]{Lemma}
\newtheorem{proposition}[theorem]{Proposition}
\theoremstyle{definition}

\theoremstyle{remark}

\numberwithin{equation}{section}
\numberwithin{theorem}{section}

\newcommand{\C}{\mathbb C}

\newcommand{\R}{\mathbb R}

\renewcommand{\H}{\mathcal H}

\def\({\left(}
\def\){\right)}

\def \mm #1#2#3#4{\begin{pmatrix} #1 & #2 \cr #3 & #4 \end{pmatrix}}

\begin{document}
\title[Converse theorems assuming a partial Euler product]
{Converse theorems assuming a partial Euler product}
\author{David W. Farmer and Kevin Wilson}

\thanks{This work resulted from an REU at Bucknell University 
and the American Institute of Mathematics.  Research supported by the
American Institute of Mathematics
and the National Science Foundation}

\thispagestyle{empty}
\vspace{.5cm}
\begin{abstract}
Associated to a newform $f(z)$ is a Dirichlet series $L_f(s)$ with
functional equation and Euler product.  Hecke showed that if the
Dirichlet series $F(s)$ has a functional equation of a particular
form, then $F(s)=L_f(s)$ for some holomorphic newform $f(z)$ on $\Gamma(1)$.
Weil extended this result to $\Gamma_0(N)$ under an assumption
on the twists of $F(s)$ by Dirichlet characters. Conrey and Farmer
extended Hecke's result for certain small~$N$, assuming that the local 
factors in the Euler product of $F(s)$ were of a special form.  
We make the same assumption on the Euler product and describe an approach 
to the converse theorem using certain additional assumptions.
Some of the assumptions may be related to second order modular forms.
\end{abstract}

\address{
{\parskip 0pt
American Institute of Mathematics\endgraf
farmer@aimath.org\endgraf
\null
University of Maryland\endgraf
kmwilson@math.umd.edu\endgraf
}
  }

\maketitle

\section{Introduction}

The connection between $L$-functions and modular
forms began with Riemann who deduced the functional equation for
the $\zeta$-function from the transformation law for the $\theta$-function.
Hecke showed that modular forms can be associated to Dirichlet series
having a  functional equation and Euler product
of a certain form.  
The converse problem, showing that $L$-functions having certain properties
arise from a modular form, is significantly more difficult.  We outline
some of the previous work on this problem and then describe our result.

The description below is somewhat imprecise, as we 
ignore certain technical issues.  See \cite{Iw} for a complete discussion,
and see Section~\ref{sec:background} for definitions
of the terms used here.

Suppose 
\begin{equation}\label{fourierexpansion}
f(z)=\sum_{n=1}^\infty a_n e(nz),
\end{equation}
where $e(z)=e^{2\pi i z}$.
We obtain a Dirichlet series
by taking the Mellin transform of $f$ along the positive imaginary axis:
\begin{eqnarray}
\xi(s,f) &=& \int_0^\infty f(i y) y^s \frac{dy}{y} \cr
&=& \(\frac{\sqrt{N}}{2\pi}\)^s \Gamma(s) L(s,f)
\end{eqnarray}
where
\begin{equation}
L(s,f)=\sum_{n=0}^\infty \frac{a_n}{n^s} .
\end{equation}
The above relationship makes it clear that there should be some
connection between properties of $f(z)$ and properties of $L(s,f)$.
Consider the following conditions:

\begin{enumerate}
\item[$A)$] $f\in S_k(\Gamma_0(N))$.

\item[$B)$]  $f$ is in the Hecke basis of $S_k(\Gamma_0(N))$.

\item[$C)$]  $f|H_N=\pm f$ and $f|T_p=a_p f$ for prime $p$. 

\item[$D)$]   $L(s,f)$ has a functional equation and an Euler product of
degree~2, weight~$k$, and level~$N$.

\item[$E)$]  For a particular set of integers $S$, the twisted $L$-function $L(s,f,\chi)$ 
has a functional equation of degree~2, weight~$k$, 
and level~$N_q$ for all characters $\bmod q$, for $q\in S$.
\end{enumerate}
 
We have the following implications:

$B\Rightarrow A$ by definition.

$A$ almost implies $B$ in the sense that the 
functions meeting condition $A$ can be written as a finite linear
combination of functions meeting condition~$B$.

$B\Rightarrow C$ by definition.

$C\Leftrightarrow D$.  This equivalence has two parts.  

\begin{itemize} 
\item The
functional equation for $L(s,f)$ is equivalent to $f|H_N=\pm f$
by using the Mellin transform and its inverse.  This observation is due
to Hecke.  (One needs $L(s,f)$ to be bounded in vertical strips in order
to justify the steps here.  See Chapter~6 of~\cite{Iw} for details).

\item The shape of the local factors in the Euler product is equivalent
to the assertion $f|T_p=a_p f$, because both
are equivalent to particular recurrence relations on the coefficients.
\end{itemize}

What is missing from the above implications is $C\Rightarrow B$.
The possibility of this was first suggested by
Conrey and Farmer~\cite{cf}, who showed that implication for
certain small~$N$.

A different approach was pioneered by Weil~\cite{w}, who showed that
$E\Rightarrow A$ provided that the set $S$ was sufficiently large 
(so in particular $D+E\Rightarrow B$).  Later versions~\cite{raz, li}
provide more general results and show more flexibility in the
choice of moduli~$S$.  Recently, Diaconu, Perelli, and Zaharescu~\cite{dpz}
have shown that $D+E\Rightarrow B$ for $S$ consisting of \emph{one}
modulus.

In this paper we consider the implication $C\Rightarrow B$, exploring
the consequences of assuming $C$ only at the prime~$p=2$.

\section{Definitions and statement of results}\label{sec:background}

We review the basic material on holomorphic modular forms.
For more information see~\cite{Iw}.

A holomorphic function $f:\H \to \C$ is a modular form 
of weight~$k$ and character~$\chi$ for $\Gamma_0(N)$ if 
$f(z)$ is regular at the cusps of $\Gamma_0(N)$ and
$$
f \left(\frac{az+b}{cz+d}\right)=\chi(d)(cz+d)^k f(z)
$$ 
for all
 $\gamma = \mm abcd \in \Gamma_0(N)$
where 
$$
\Gamma_0(N)=\left\{\mm abcd : a,b,c,d \hbox{ are integers, } 
ad - bc = 1 \hbox{, and } c \equiv 0 \bmod N \right\}
$$
and $\H= \{z =x+iy\in \C : y>0\}$ is the upper half-plane.
Furthermore, $f$ is a cusp form if also it 
vanishes at all cusps of $\Gamma_0(N)$.  The collection of cusp 
forms of weight $k$ and character $\chi$ for $\Gamma_0(N)$
forms a vector space which is denoted $S_k(\Gamma_0(N),\chi)$.  We set 
$S_k(\Gamma_0(N))=S_k(\Gamma_0(N),\chi_0)$ where
$\chi_0$ is the trivial character mod~$N$.

We will now introduce some notation.  Let $\gamma=\mm abcd \in SL(2,\R)$.
The following are four equivalent ways of stating the same information:
\begin{eqnarray}
f \left(\frac{az+b}{cz+d}\right) &=& \chi(d)(cz+d)^{k} f(z) \cr
f( \gamma z) &=& \chi(d)(cz+d)^{k} f(z) \cr
f|\gamma&=&f \cr
\gamma &\equiv& 1 \bmod \Omega_f.
\end{eqnarray}
Here $\Omega_f$ is the right ideal in the group ring $\C[GL(2,\R)]$
which annihilates~$f$, the action of matrices on $f$ being extended
linearly.
We will write $\gamma \equiv 1$ instead of $\gamma \equiv 1 \bmod \Omega_f$ 
from Section~\ref{proofofyougetG02} onward.  
In the third displayed line above we used the standard ``slash'' notation
\begin{equation}
(f|\gamma)(z) = (\det\gamma)^{-k/2} \chi(d)(cz+d)^{k} f(z),
\end{equation}
which will be used below as we make the discussion from the introduction
more precise.


The vector space $S_k(\Gamma_0(N))$ has a distinguished basis of
elements, called the Hecke basis, which satisfy
\begin{equation}\label{frickepm}
f|H_N=\pm f
\end{equation}
and
\begin{equation}
f|T_p=a_p f
\end{equation}
for prime $p$.
Here
$$
H_N = \mm {}{-1}N{} 
$$
is the Fricke involution.  If $\ell$ is prime,
$$
T_\ell = \chi(\ell) \mm \ell{}{}1 + \sum_{a=0}^{\ell-1} \mm 1a{}\ell,
$$ 
is the Hecke operator.  If $\ell$ divides $N$ then $\chi(\ell)=0$ and $T_\ell$ is known
as the Atkin-Lehner operator~$U_\ell$. 

An $L$-function is a Dirichlet Series with a functional equation and 
an Euler product.  
In the case that $L(s)$ is associated to a newform in  $S_k(\Gamma_0(N))$
then we have the functional equation
\begin{equation} \label{thefunctionalequation}
\xi_f(s)=\(\frac{2\pi}{\sqrt{N}}\)^{-s}\Gamma\(s+\frac{k-1}{2}\)L_f(s)=\pm (-1)^{k/2}\xi_f(1-s).
\end{equation}
and an Euler product of the form
\begin{equation} \label{theEulerproduct}
L(s,f)=\prod_p
\(1-a_p p^{-s} + \chi(p) p^{k-1-2s}\)^{-1}
\end{equation}

The purpose of this paper is to determine what information could possibly 
be deduced from the functional equation (\ref{thefunctionalequation}) and
an Euler factor only at the prime~2:
\begin{equation} \label{partialEulerproduct}
L(s,f)=\(1-a_2 2^{-s} + 2^{k-1-2s}\)^{-1}
\sum_{n\ odd} \frac{a_n}{n^s}
\end{equation}
The functional equation (\ref{thefunctionalequation}) is equivalent
to
$f|H_N = \pm f$.
On the other hand, the Euler product (\ref{partialEulerproduct})
is equivalent to 
$f|T_2=a_2f$, where $\chi(2)=1$.   So it may be that 
$f\in S_k(\Gamma_0(N))$ for some odd~$N$, or possibly that
$f\in S_k(\Gamma_0(N),\chi)$ where $\chi(2)=1$.
But more generally it may only be that $f\in S_k(\Gamma_{0,2}(N))$,
where
$$
\Gamma_{0,2}(N)=\left \{\mm abcd \in \Gamma_0(N): 
a \equiv 2^j \bmod N \hbox{ for some integer j} \right \}.
$$
These ambiguities are unavoidable and must factor into our results.

We will make repeated use of the fact that if $f$ has a Fourier expansion
of the form~(\ref{fourierexpansion}) then
$f(z)=f(z+1)$, and if $f$ also satisfies~(\ref{frickepm}) then 
$f|T=f$ and $f|W_N=f$ where
$$
T=\mm 11{}1 \hbox{ and } W_N= H_N T^{-1} H_N^{-1} = \mm 1{}N1.
$$

Since $T$ and $W_N$ generate $\Gamma_0(N)$ for $N\le 4$,
in the notation of the introduction we have $C\Rightarrow B$ for~$N\le 4$,
which is Hecke's converse theorem.  Note that this does
not make use of the Euler product.
Our goal is to examine the implication
$C\Rightarrow B$ for larger $N$ under the additional condition
$f|T_2=a_2f$.  

The following is our main result.  It is conditional on certain
assumptions which are described later in the paper.

\begin{theorem}\label{maintheorem}
Under Pairing Assumptions 1 and 2 (described in Section~\ref{M2case}), 
and a weak form of Artin's conjecture on primitive roots (described in
Section~\ref{proofthatG02generates}),
if $f|H_N=\pm f$, $f|T=f$, and $f|T_2=a_2f$, then $f|\gamma=f$ for all
$\gamma$ in the group
$$
\Gamma_{0,2}(N)=\left \{\mm abcd \in \Gamma_0(N): 
a \equiv 2^j \bmod N \hbox{ for some integer j} \right \}.
$$
\end{theorem}

Note that 
$$
\Gamma_{0,2}(N) = \left\{\gamma \in \Gamma_0(N): \chi(\gamma)=1 
\hbox{ for all characters mod $N$ such that } \chi(2)=1 \right\}.
$$
By the comments following (\ref{partialEulerproduct})
we see that 
$\Gamma_{0,2}(N)$ is in fact the largest group for which one could
hope to conclude $f|\gamma=f$ from our assumptions.
Note also that $\Gamma_{0,2}(N)$ contains $\Gamma_{1}(N)$, in particular
it has finite index in~$\Gamma(1)$.

The proof of Theorem~\ref{maintheorem} will involve two main steps.  
First, we show

\begin{proposition}\label{yougetG02} Under Pairing Assumptions 1 and 2, 
if $f|H_N= \pm f$, $f|T=f$, and $f|T_2=a_2f$
then $f|\gamma=f$ for all $\gamma$ in the set 
$$
G_{0,2}(N)= \left\{\mm {2^n}\alpha{\beta N}{*} \in \Gamma_0(N)\ :\ n\ge0\right\}.
$$
\end{proposition}
Here and following, a single entry of $*$ in a matrix indicates that the matrix
has determinant~1, which can be used to determine the value of~$*$.

Once we have Proposition~\ref{yougetG02}, then Theorem~\ref{maintheorem} follows from 

\begin{proposition}\label{G02generates} Assuming a weak form of Artin's conjecture
on primitive roots,
$G_{0,2}(N)$ is a generating set for $\Gamma_{0,2}(N)$. 
\end{proposition}

We believe that the version of Artin's conjecture required for the Proposition
is accessible by current methods, and furthermore it should be possible to
prove the proposition without that assumption.

We prove Proposition~\ref{yougetG02} in the next section, and we prove 
Proposition~\ref{G02generates} in Section~\ref{proofthatG02generates}.

\section{Proof of Proposition~\ref{yougetG02}}\label{proofofyougetG02}

The first step is to combine the information $H_N \equiv \pm 1$, 
$T \equiv 1$,  
and $T_2 \equiv a_2$ to show that $\gamma \equiv 1$ 
for all $\gamma$ in the set
$$
L_{0,2}(N)=\left\{\mm {2^n}\alpha{\beta N}{*} \in \Gamma_0(N):|\alpha|,|\beta| \le 2^{n-1}\right\}.
$$
We will then be done with the majority of the proof of Proposition~\ref{yougetG02} 
and be left with the easy task of showing that $T$, $W_N$, 
and $L_{0,2}(N)$ generate the set $G_{0,2}(N)$.

\subsection{Base case}

We will now show that for $\gamma = 
\mm {2^n}\alpha{\beta N}{*} \in L_{0,2}(N)$, 
if $n=1$ then $\gamma \equiv 1$.

\begin{lemma}[Lemma 2 of \cite{cf}]  If $H_N \equiv \pm 1$ and 
$T_2\equiv a_2$ then 
$$
\mm 2{\alpha}{\beta N}{\frac{\alpha \beta N+1}{2}} 
\equiv 1
$$ 
if $|\alpha|=|\beta|=1$.
\end{lemma}

\begin{proof}
Assume $H_N \equiv \pm 1$ and $T_2\equiv a_2$.  
Note that:
$$
H_N^{-1} T_2 H_N 
=  \mm 1{}{}2 + \mm 2{}{}1 + \mm 2{}{-N}1.
$$

Since $H_N^{-1} T_2 H_N \equiv a_2 H_N^{-1} H_N \equiv a_2\equiv T_2$, we have:
$$
\mm 1{}{}2 + \mm 2{}{}1 + \mm 2{}{-N}1 \equiv 
\mm 2{}{}1 + \mm 1{}{}2 + \mm 11{}2.
$$

Canceling common terms from both sides we are left with:
$$
\mm 2{}{-N}1 \equiv \mm 11{}2.
$$
Right multiplying by $\mm 11{}2^{-1}$ we have:
$$
M_2:=\mm 2{-1}{-N}{\frac{N+1}{2}} \equiv 1.
$$
The other matrices in $L_{0,2}(N)$ with $n=|\alpha|=|\beta|=1$ are 
obtained from $T$, $W_N$, and $M_2$ in the following way:
\begin{eqnarray}
\mm 2{-1}N{\frac{-N+1}{2}} &=& W_N \mm 2{-1}{-N}{\frac{N+1}{2}},\cr
\mm 21{-N}{\frac{-N+1}{2}} &=& \mm 2{-1}{-N}{\frac{N+1}{2}} T,\cr
\end{eqnarray}
and
$$
\mm 21N{\frac{N+1}{2}} = W_N \mm 2{-1}{-N}{\frac{N+1}{2}} T.
$$
Thus we have that all of these different forms are equivalent to 1, 
as required.  
\end{proof}

Now that we have $M_2 \equiv 1$, we will use this to produce
more matrices $\gamma \equiv 1$.  This is analogous to the 
above process where $H_N\equiv 1$ led to $M_2 \equiv 1$.

\subsection{First induction step}\label{M2case}

We now start with
$$
M_2 = \mm 2{-1}{-N}{\frac{N+1}{2}} \equiv 1,
$$ 
which implies
$$
T_2 M_2 \equiv a_2 M_2 \equiv a_2 \equiv T_2.
$$ 
Note that:
$$
T_2 M_2 
= 
\mm 4{-2}{-N}{\frac{N+1}{2}} +
\mm 2{-1}{-2N}{N+1} +
\mm {2-N}{\frac{N-1}{2}}{-2N}{N+1}.
$$ 
So by $T_2 M_2 \equiv T_2$ we have:
$$
\mm 4{-2}{-N}{\frac{N+1}{2}} +
\mm 2{-1}{-2N}{N+1} +
\mm {2-N}{\frac{N-1}{2}}{-2N}{N+1}
\equiv
\mm 2{}{}1 + \mm 1{}{}2 + \mm 11{}2.
$$
We would like four of these matrices to cancel, 
as occurred in the $n=1$ case. 

Recall that $T \equiv 1$
and $W_N \equiv 1$.  Notice that 
$$
W_N^{-1} T^{-1} \mm 2{}{}1 =  \mm 2{-1}{-2N}{N+1} ,
$$
so we have 
$$
\mm 2{-1}{-2N}{N+1} \equiv \mm 2{}{}1.
$$
Thus, two of those matrices cancel,  leaving
\begin{equation}\label{fourmatrices}
\mm 4{-2}{-N}{\frac{N+1}{2}}+\mm {2-N}{\frac{N-1}{2}}{-2N}{N+1} \equiv
\mm 1{}{}2+\mm 11{}2.
\end{equation}
Unfortunately, for general $N$ we cannot cancel any more matrices, 
so things do not work out quite as nicely as they did when using $H_N$.  
This leaves us at an impasse so we need to make an assumption. 
We will now state our first pairing assumption.

\medskip
\noindent{\bf Pairing Assumption 1. } \emph{If 
$A+B \equiv C+D$, then either $A \equiv C$ and $B \equiv D$, or 
$A \equiv D$ and $B \equiv C$.}

\medskip

On Pairing Assumption 1, we have that the
matrices in (\ref{fourmatrices}) must pair up one of two ways.  
Either 
$$
\mm 4{-2}{-N}{\frac{N+1}{2}} \equiv \mm 1{}{}2 \hbox{ and } 
\mm {2-N}{\frac{N-1}{2}}{-2N}{N+1} \equiv \mm 11{}2
$$
which gives 
$$
\mm 4{-1}{-N}{\frac{N+1}{4}} \equiv 1 \hbox{ and } 
\mm 41N{\frac{N+1}{4}} \equiv 1
$$
or
$$
\mm 4{-2}{-N}{\frac{N+1}{2}} \equiv \mm 11{}2 \hbox{ and } 
\mm {2-N}{\frac{N-1}{2}}{-2N}{N+1} \equiv \mm 1{}{}2
$$
which gives 
$$
\mm 4{-1}N{\frac{-N+1}{4}} \equiv 1 \hbox{ and } 
\mm 41{-N}{\frac{-N+1}{4}} \equiv 1.
$$

If $N \equiv 3 \bmod 4$ then the first way of 
pairing  gives integer matrices and 
if $N \equiv 1 \bmod 4$ then the second way gives 
integer matrices.
As we will see in Lemma~\ref{oneisinteger},  
in general one of the pairings gives integer matrices and the other does not.  
As a result of this, we will need to make another assumption.

\medskip
\noindent{\bf Pairing Assumption 2. } \emph{ If 
$A+B \equiv C+D$,
and Pairing Assumption 1 is true, then $A,B,C,D$ pair in 
a way that gives integer matrices equivalent to~1. }

\medskip
So we have now gone through the process using $T_2 \equiv a_2$ and 
$M_2 \equiv 1$ which, by the Pairing Assumptions,  gave two
more matrices in $L_{0,2}(N)$ which are equivalent to 1.  
We will now generalize the sequence of steps we have just used in the 
cases $n=1$ and~$2$.

\subsection{The Process}

In the previous section we used $T_2M_2\equiv T_2$ to obtain, 
under the pairing assumptions, $\gamma\equiv1$ for two
new matrices~$\gamma$.
We can now try to use $T_2\gamma\equiv T_2$, that is 
\begin{equation}\label{T2gammaT2}
\mm 2{}{}1\gamma +\mm 1{}{}2\gamma +\mm 11{}2\gamma \equiv 
\mm 2{}{}1 + \mm 1{}{}2 + \mm 11{}2
\end{equation}
to obtain
two more matrices $\gamma'\equiv 1$.  
If this step is successful then we can continue the process.
Our goal is to understand what matrices can be obtained in this way,
assuming Pairing Assumptions~1 and~2. 

\subsection{Producing Matrices from the Process}

We will show that all matrices produced from the process 
are in $L_{0,2}(N)$ and that we actually obtain all of the matrices in 
$L_{0,2}(N)$.

\begin{proposition}\label{twomoreinL2} Under Pairing Assumptions 1 and 2,
if $\gamma\equiv1$ for $\gamma \in L_{0,2}(N)$ with $N$ odd, 
and you repeat the process
$T_2\gamma\equiv T_2$,
then you obtain two more matrices $\gamma_1, \gamma_2 \in L_{0,2}(N)$, 
such that $\gamma_1 \equiv \gamma_2 \equiv 1$.
\end{proposition}

\begin{proof} We will show by induction  that the process 
leads to $\gamma \equiv 1$ for all $\gamma = \mm{2^n}{\alpha}{\beta}{*} \in L_{0,2}(N)$.  

The base case $n=1$ comes from
$T_2 H_N \equiv T_2$ which yielded $M_2 \equiv 1$.  
 
Now suppose $\mm {2^j} \alpha {\beta N} {*}
\equiv 1$ for all $j \le n$ and for all $|\alpha|$, $|\beta| \le 2^{j-1}$.  
Let $\gamma = \mm {2^{n}}{\alpha}{\beta N}{*}$ 
where $\alpha$ and $\beta$ are odd and $|\alpha|,|\beta| \le 2^{n-1}$.
We have $T_2\gamma\equiv T_2$, so
\begin{eqnarray}
T_2 \gamma &=& \mm 2{}{}1 \gamma + \mm 1{}{}2 \gamma
 + \mm 11{}2 \gamma \cr 
&= &
\mm {2^{n+1}}{2\alpha}{\beta N}{\frac{\alpha \beta N+1}{2^n}}+
\mm {2^{n}}{\alpha}{2\beta N}{\frac{\alpha \beta N+1}{2^{n-1}}}+
\mm {2^{n}+\beta N}{\alpha+\frac{\alpha \beta N+1}{2^{n-1}}}
{2\beta N}{\frac{\alpha \beta N+1}{2^{n-1}}} \cr 
&\equiv &
\mm 2{}{}1 + \mm 1{}{}2 + \mm 11{}2  .
\end{eqnarray}
We expect that one pair of matrices will cancel, as it did
in the previous example.

By the induction hypothesis $\mm {2^{n-1}}\alpha{\beta N}
{\frac{\alpha \beta N+1}{2^{n-1}}} \equiv 1$.
Also note that
$$
\mm {2^{n}}{\alpha}{2\beta N}{\frac{\alpha \beta N+1}{2^{n-1}}}
=
\mm {2^{n-1}}\alpha{\beta N}{\frac{\alpha \beta N+1}{2^{n-1}}}
\mm 2{}{}1
$$ 
so we have 
$$
\mm {2^{n}}{\alpha}{2\beta N}
{\frac{\alpha \beta N+1}{2^{n-1}}} \equiv \mm 2{}{}1
$$
so those matrices cancel in the previous equivalence.
This leaves four matrices:
$$
\mm {2^{n+1}}{2\alpha}{\beta N}{\frac{\alpha \beta N+1}{2^n}}+
\mm {2^{n}+\beta N}{\alpha+\frac{\alpha \beta N+1}{2^{n-1}}}
{2\beta N}{\frac{\alpha \beta N+1}{2^{n-1}}}  \equiv
\mm 1{}{}2 + \mm 11{}2 .
$$
There are now two possible ways for the matrices 
to pair up.  So by Pairing Assumption~1 we have either
\begin{eqnarray}
\mm {2^{n+1}}{2\alpha}{\beta N}{\frac{\alpha \beta N+1}{2^n}}
\mm 1{}{}2^{-1}&=&
\mm {2^{n+1}}{\alpha}{\beta N}{\frac{\alpha \beta N+1}{2^{n+1}}} \cr
&\equiv & 1
\end{eqnarray}
and
\begin{eqnarray}
\mm {2^{n}+\beta N}{\alpha+\frac{\alpha \beta N+1}
{2^{n-1}}}{2\beta N}{\frac{\alpha \beta N+1}{2^{n-1}}}
\mm 11{}2^{-1}&=&
\mm {2^{n}+\beta N}{\frac{\alpha \beta N-2^n\beta N-2^{2n}+2^n\alpha+1}
{2^{n+1}}}{2\beta N}{\frac{\alpha \beta N-2^n \beta N+1}{2^{n}}} \cr
&\equiv & M_2 \mm {2^{n}+\beta N}
{\frac{\alpha \beta N-2^n\beta N-2^{2n}+2^n\alpha+1}
{2^{n+1}}}{2\beta N}{\frac{\alpha \beta N-2^n \beta N+1}{2^{n}}} \cr
&=& \mm {2^{n+1}}{\alpha-2^n}{(\beta-2^n) N}{\frac{(\alpha-2^n)(\beta-2^n)N+1}
{2^{n+1}}}  \cr
&\equiv & 1,
\end{eqnarray}
or
\begin{eqnarray}
\mm {2^{n+1}}{2\alpha}{\beta N}{\frac{\alpha \beta N+1}{2^n}}
\mm 11{}2^{-1}&=&
\mm {2^{n+1}}{\alpha-2^n}{\beta N}{\frac{(\alpha-2^n)\beta N+1}
{2^{n+1}}} \cr
&\equiv & 1 
\end{eqnarray}
and
\begin{eqnarray}
\mm {2^{n}+\beta N}{\alpha+\frac{\alpha \beta N+1}
{2^{n-1}}}{2\beta N}{\frac{\alpha \beta N+1}{2^{n-1}}}
\mm 1{}{}2^{-1}&=&
\mm {2^{n}+\beta N}{\frac{\alpha \beta N+2^n\alpha+1}{2^{n+1}}}
{2\beta N}{\frac{\alpha \beta N+1}{2^{n}}} \cr
&\equiv & M_2 \mm {2^{n}+\beta N}{\frac{\alpha \beta N+2^n\alpha+1}{2^{n+1}}}
{2\beta N}{\frac{\alpha \beta N+1}{2^{n}}} \cr
&=& \mm {2^{n+1}}{\alpha}{(\beta - 2^n)N}
{\frac{\alpha(\beta - 2^n) N+1}{2^{n+1}}} \cr
&\equiv & 1.
\end{eqnarray}

There are now two things left to be shown.  First, we must show 
that the matrices we get from this process are actually of the 
form of matrices in $L_{0,2}(N)$.  That is, we must show that the
off-diagonal entries can be made small enough.   Then we must show that 
one of these pairs are integer matrices.  

\begin{lemma}\label{niceform} With $\alpha$, $\beta$, and $n$ defined as above,
${\mm {2^{n+1}}{\alpha-2^n}{(\beta-2^n) N}
{*}}$,  
${\mm {2^{n+1}}{\alpha-2^n}{\beta N} 
{*}}$, 
and
${\mm {2^{n+1}}{\alpha}
{(\beta-2^n) N}{*}}$ 
can be put in the form
$\mm {2^{n+1}}{\alpha'}{\beta' N}
{*}$ with $|\alpha'|,|\beta'|<2^n$,
using only multiplication by $T$ and $W_N$.
\end{lemma}

\begin{proof} 
Note that 
$$
\mm ab{cN}* T^\ell = \mm a{b+a\ell}{cN}{*} \hbox{ and } W_N^k \mm ab{cN}* = 
\mm ab{(c+ak)N}{*}.
$$
Now choose $\ell$ and $k$ 
such that  
$|b+2^{n+1} \ell|$, $|c + 2^{n+1}k| \le 2^n$.  
Thus we can always adjust the value of $b$ and $c$ 
by $2^{n+1}=2*2^n$ without changing the top left entry.  
This is sufficient because for any 
$|x| < 2^{n+1}$, one of $|x|$, $|x-2^{n+1}|$, 
or $|x+2^{n+1}|$ is less than~$2^n$.  
\end{proof}

All that is left to be shown is that one pair of matrices will
have integer entries.
From Pairing Assumption~2, we pair them in the 
way that gives integer matrices equivalent to one.
The following lemma shows that this happens for exactly
one of the pairings.

\begin{lemma}\label{oneisinteger}  If $N$ is odd and 
$\mm {2^{n}}\alpha{\beta N}{\frac{\alpha \beta N+1}{2^{n}}} \in \Gamma_{0}(N)$, 
then exactly one of the pairs
$$
\mm {2^{n+1}}\alpha{\beta N}{\frac{\alpha \beta N+1}{2^{n+1}}}
\ \ \ and\ \ \  
\mm {2^{n+1}}{\alpha - 2^n}{(\beta - 2^n)N}
{\frac{(\alpha - 2^n)(\beta - 2^n) N+1}{2^{n+1}}}
$$
or
$$
\mm {2^{n+1}}{\alpha - 2^n}{\beta N}
{\frac{(\alpha - 2^n) \beta N+1}{2^{n+1}}}
\ \ \ and\ \ \ 
\mm {2^{n+1}}{\alpha}{(\beta - 2^n)N}
{\frac{\alpha(\beta - 2^n) N+1}{2^{n+1}}}
$$
is in $\Gamma_{0}(N)$.
\end{lemma}

\begin{proof}
The issue is whether the bottom right entry of the matrices is an integer.  
It is straightforward to check that if
$\frac{\alpha \beta N+1}{2^n}$ 
is even then the first pair are in $\Gamma_0(N)$, otherwise
the second pair is.
\end{proof}
 
Since we showed in Lemma~\ref{niceform} that 
the entries in these matrices 
can be made of the appropriate size, we have that they 
are actually in $L_{0,2}(N)$, which completes the proof of 
Proposition~\ref{twomoreinL2}.
\end{proof}

We have shown that the process
leads to matrices in $L_{0,2}(N)$.  
It remains to show that all of the matrices in 
$L_{0,2}(N)$ can be obtained.

\begin{proposition} If $\alpha$ and $\beta$ are odd and 
$|\alpha|,|\beta| \le 2^{n-1}$, then there exists a matrix 
$\gamma=\mm {2^{n-1}}{\alpha_1}{\beta_1 N}
{\frac{\alpha_1 \beta_1 N + 1}{2^{n-1}}} \in L_{0,2}(N)$ 
such that if $\gamma\equiv 1$ and $T_2\gamma\equiv T_2$,
then Pairing Assumptions~1 and~2 lead to
$\mm {2^{n}}\alpha{\beta N}
{\frac{\alpha \beta N + 1}{2^{n}}}\equiv 1$.
\end{proposition}

Thus, under our pairing assumptions we have $\gamma \equiv 1$ 
for all $\gamma \in L_{0,2}(N)$.

\begin {proof}  Let $A=\mm {2^{n}}\alpha{\beta N} {*}$. 
Now choose $\alpha_1, \beta_1 \le 2^{n-2}$ such that 
$\alpha \equiv \alpha_1 \bmod 2^{n-1}$ and 
$\beta \equiv \beta_1 \bmod 2^{n-1}$.  
From the proof of Proposition~\ref{twomoreinL2}
starting with
$\mm {2^{n-1}}{\alpha_1}{\beta_1 N}
{\frac{\alpha_1 \beta_1 N + 1}{2^{n-1}}} \equiv 1$ we obtain the
matrices:
\begin{eqnarray}
\gamma_1&= & \mm {2^{n}}{\alpha_1}{\beta_1 N}
{\frac{\alpha_1 \beta_1 N+1}{2^{n}}} \cr  
\gamma_2&= & \mm {2^{n}}{\alpha_1}{(\beta_1 - 2^{n-1})N}
{\frac{\alpha_1(\beta_1 - 2^{n-1}) N+1}{2^{n}}} \cr
\gamma_3&= & \mm {2^{n}}{\alpha_1 - 2^{n-1}}{\beta_1 N}
{\frac{(\alpha_1 - 2^{n-1}) \beta_1 N+1}{2^{n}}} \cr 
\gamma_4&= & \mm {2^{n}}{\alpha_1 - 2^{n-1}}{(\beta_1 - 2^{n-1})N} 
{\frac{(\alpha_1 - 2^{n-1})(\beta_1 - 2^{n-1}) N+1}{2^{n}}}.
\end{eqnarray}
We must show that $A\equiv \gamma_j$ for some $j$.
This breaks down into cases.

If $|\alpha|,|\beta|<2^{n-2}$, then we can set $\alpha=\alpha_1$ and 
$\beta=\beta_1$, so $A=\gamma_1$. 

Suppose $|\alpha|<2^{n-2}$. If $2^{n-2}<\beta<2^{n-1}$, then
$A=W_N \gamma_2$, and if $-2^{n-1}<\beta<2^{n-2}$ then
$A=\gamma_2$.

Suppose $|\beta|<2^{n-2}$. If $2^{n-2}<\alpha<2^{n-1}$, then
$A=\gamma_3 T$, and if $-2^{n-1}<\alpha<2^{n-2}$ then
$A=\gamma_3$.

In the four cases where $2^{n-2}<|\alpha|,|\beta|<2^{n-1}$,
use $\gamma_4$, $T$, and $W_N$ in a manner similar to the above.

Since we have exhausted all possibilities of $\alpha$ and $\beta$, 
we are done.
\end{proof}

Now we have that $L_{0,2}(N)$ is the set of matrices 
generated by the process so we are finished with the first part 
of Proposition~\ref{yougetG02}.  Now we just need to show that $T$, $W_N$, and 
$L_{0,2}(N)$ generate $G_{0,2}(N)$.  

\begin{proposition} $T$, $W_N$, and 
$L_{0,2}(N)$ generate $G_{0,2}(N)$.
\end{proposition}

\begin{proof}  We will show 
that an arbitrary matrix in $G_{0,2}(N)$ can be written in terms of 
the matrices in $L_{0,2}(N)$ along with $T$ and $W_N$.  Let 
$\psi \in G_{0,2}(N)$, so $\psi = \mm {2^{n}}{\alpha}{\beta N}
{\frac{\alpha \beta N+1}{2^{n}}}$ for some odd $\alpha$ and $\beta$.  
Now choose $\ell$ and $k$ such that $|\alpha +2^n \ell|<2^{n-1}$ and 
$|\beta + 2^nk|<2^{n-1}$.  Then 
$$
W_N^k \psi T^\ell= \mm {2^{n}}{\alpha_1}{\beta_1 N}
{\frac{\alpha_1 \beta_1 N+1}{2^{n}}} \in L_{0,2}(N),
$$
because $|\alpha_1|$, $|\beta_1| \le 2^{n-1}$.  
\end{proof}

So we now have that, starting with $T_2 \equiv a_2$, $H_N \equiv 1$, 
and $T \equiv 1$, the Pairing Assumptions give $\gamma \equiv 1$ for all 
$\gamma \in G_{0,2}(N)$. This finishes the proof of Proposition~\ref{yougetG02}.

It remains to show that $G_{0,2}(N)$ generates $\Gamma_{0,2}(N)$.

\section{Generating $\Gamma_{0,2}(N)$}\label{proofthatG02generates}

In order to prove Proposition~\ref{G02generates}, 
we will need to show that 
we can write every matrix in $\Gamma_{0,2}(N)$ 
in terms of matrices in $G_{0,2}(N)$.  
Our proof uses Artin's Conjecture in the following form:

\medskip
\noindent{\bf A weak form of Artin's Conjecture on primitive roots. } 
\emph{If $(d,bM)=1$ then there exists integers $k$ and $n$ such that
$b\equiv 2^n\mod (d+k M)$.}
\medskip

Note that this follows from Artin's conjecture, for if $p=d+k M$
was prime and $2$ was a primitive root mod~$p$, then the conclusion
would follow.  But the above statement is actually much weaker, for
there is no requirement that $2$ be a primitive root, nor that
$p$ be prime.  Indeed, on average the powers  of $2$ are
a substantial fraction of the integers mod~$d+k M$, so barring
any obstruction one would expect the subgroup generated by 2 to
contain $b$ fairly often.

\begin{proof}[Proof of Proposition~\ref{G02generates}]  
Let $\gamma= \mm abcd \in \Gamma_{0,2}(N)$,  and let $a \equiv 2^j \bmod N$.
First suppose $b\not=1$ and
let $\delta=\mm{2^n}{\alpha}{*}{*} \in G_{0,2}(N)$, so
we have 
$$
\delta W_N^{k} \gamma = \mm *{2^n b+ \alpha(d+b N k)} * * .
$$
By the weak form of Artin's Conjecture, choose $k$ and $n$ such 
that $2^n b+ \alpha(d+b N k)=1$.
Thus, we have reduced to the case $\gamma= \mm abcd \in \Gamma_{0,2}(N)$
with $b=1$. 

Now let $\gamma = \mm a1{ad-1}d \in \Gamma_{0,2}(N)$ where 
$a \equiv 2^j \bmod N$.  
Since $W_N \in G_{0,2}(N)$, it is sufficient to show that 
$\gamma W_N^{k} \in G_{0,2}(N)$ for some $k$.  

To have $\gamma W_N^{k} \in G_{0,2}(N)$, we just need 
the top left entry to be a power of two.  
Since $a \equiv 2^j \bmod N$, we have $a+k N=2^j$ for 
some integer $k$.  So 
$$
\gamma W_N^{k} = \mm {2^j} 1 {ad-1 + d k N} d \in G_{0,2}(N),
$$
as required.
\end{proof}

This completes the proof of Theorem~\ref{maintheorem}, assuming a weak
form of Artin's Conjecture.
Note that the proof actually shows that every element of $\Gamma_{0,2}(N)$
is of the form $W_N^k \delta_1 \delta_2 W_N^\ell$, where
$\delta_1,\delta_2 \in G_{0,2}(N)$.
By considering more general products of elements in $ G_{0,2}(N)$, it should
be possible to avoid the use of Artin's conjecture.

\section{Comments on the Pairing Assumptions}

Theorem~\ref{maintheorem}  reduces the 
converse theorem to showing the Pairing Assumptions. 
It is likely that the Pairing Assumptions are not true in complete
generality.  However, when we apply them in this paper, it is in the
context of functions which satisfy $f|T=f$, $f|H_N=\pm f$, and
$f|T_2=a_2 f$.  Thus, it is possible that there are extra conditions
on $f$ which make the pairing assumptions true,  and these could be deduced
from the three original assumptions.

There are situations in which 
Pairing Assumption~1 can be shown to hold.  For example, in the proof
of Weil's converse theorem (this also is used in the paper of Conrey and Farmer),
a key step is the lemma that if $\varepsilon$ is an elliptic matrix
of infinite order, then $f|(1-\gamma)(1-\varepsilon)=0$ implies
$f|(1-\gamma)=0$.  In the setup of Pairing Assumption~1, this
means that $1 + \gamma \varepsilon \equiv \gamma + \varepsilon$
implies $1 \equiv \gamma$, which is the conclusion we seek.
The examples in this paper do not appear to reduce to that situation,
but Weil's lemma can be viewed as support for the possibility
that something like Pairing Assumption~1 is true.

We believe that showing Pairing Assumption 2 is less difficult than 
showing Pairing Assumption~1.  
We have actually made some progress on showing Pairing Assumption~2 
but we have not been able to finish the proof.  
The argument begins with the observation that if you make 
a choice that gives non-integer matrices,
then as you continue 
the process you continue to generate more non-integer matrices.
The plan is to show that you eventually generate a group which 
is not discrete.
This would be a contradiction because we assumed that $f$ was a nonconstant
holomorphic function.
The issue is to
prove that the infinite set of non-integer matrices generates
a non-discrete group.

We believe Pairing Assumption 1 to be a much harder problem  
whose solution cannot be a purely algebraic 
one.  The solution cannot be algebraic because the way that you 
pair the matrices to get integers entries depends on~$N$.  
It is possible that a counterexample to Pairing Assumption 1 
could be found among second order modular forms~\cite{CDO,DKMO}.
These are functions which satisfy $f|(1-\gamma)(1-\delta)=0$ for 
$\delta \in \Gamma$.  
This 
implies $f|1+\gamma \delta -\gamma -\delta=0$ which then gives 
$1+\gamma \delta \equiv \gamma + \delta$.
If $f|(1-\gamma)\not\equiv 0$ then the terms in this equation
may not pair up.  However, it is not clear if second order modular
forms are able to satisfy the other assumptions we put on
our functions.

\end{document}